\documentclass[12pt]{article}
\usepackage[cp1251]{inputenc}
\usepackage[russian]{babel}
\usepackage{amssymb,latexsym,amsmath}
\usepackage{array,longtable}
\usepackage{enumerate}

\newcommand{\D}{\displaystyle}
\newcommand{\dok}{{\bf Proof.}\ }

\newcounter{th}

\newcommand{\theorem}{\par\refstepcounter{th}%
\bf Theorem~\arabic{th}\ }

\newcounter{lm}

\newcommand{\lem}{\par\refstepcounter{lm}%
\bf Lemma~\arabic{lm}\ }
\newenvironment{lemma}{\lem\it{}{}}

\makeatletter

\renewcommand{\section}{\@startsection{section}{1}{0pt}%
  {-3.5ex plus -1ex minus -.2ex} {2.3ex plus .2ex}%
  {\normalfont\large\bfseries}}
\renewcommand{\subsection}{\@startsection{subsection}{2}{0pt}%
  {-3.5ex plus -1ex minus -.2ex} {1.5ex plus .2ex}%
  {\normalfont\bfseries}}

\begin{document}
\sloppy

\begin{center}
{\large\bf THE GENERALIZED NEWTON--KANTOROVICH METHOD FOR
EQUATIONS WITH NONDIFFERENTIABLE OPERATORS}
\end{center}

\begin{center}
{\large\bf A.~N.~Tanyhina}
\end{center}

\begin{center}
{\it Department of Mechanics and Mathematics, \\
Belarusian State University, \\
220030 Minsk, Belarus \\
{\rm e-mail: anast-minsk@yandex.ru}}
\end{center}

\bigskip

{\small

\textbf{Abstract:} The article deals with the generalized
Newton--Kantorovich method for solving operator equations with
nondifferentiable operators in Banach spaces. The convergence
theorem is proved by means of majorant scalar equations.

\medskip


\medskip


}

\section{Introduction}

Let $X$ and $Y$ be Banach spaces, $D$ is a convex subset of $X,$ $f$
and $g$ are nonlinear operators, defined on $D$ and taking values in
$Y,$ where $f$ is differentiable at every interior point of $D,$ $g$
is nondifferentiable. One of the most effective iterative method for
solving operator equation of the form
\begin{equation}\label{osnur}
f(x) + g(x) = 0
\end{equation}
is the generalized Newton--Kantorovich method with successive
approximations
\begin{equation}\label{priblNewton}
x_{n+1} = x_n - [f'(x_n)]^{-1}(f(x_n) + g(x_n)) \quad (n = 0,\ 1,\
\ldots),
\end{equation}
where $x_0 \in D$ is given.

A thorough convergence analysis of the sequence
\eqref{priblNewton} was carried in \cite{ZabrZlep} by means of the
approach based on the application of majorant scalar equations and
originating from Kantorovich's investigations (\cite{Kantorovich},
chapter~XVIII). However the hypotheses given there are tediously
formulated and difficult to verify. For this reason in \cite{ZF}
was proposed a more flexible approach for solving the equation
\eqref{osnur} under the following hypotheses on the operators $f$
and $g$:
\begin{equation}\label{u1}
\|f'(x'') - f'(x')\| \leq \varphi(t)\|x'' - x'\|,\quad \forall\,x',\
x'' \in \overline{B(x_0, t)} \subseteq D,
\end{equation}
\begin{equation}\label{u2}
\|g(x'') - g(x')\| \leq \psi(t)\|x'' - x'\|,\quad \forall\,x',\
x'' \in \overline{B(x_0, t)} \subseteq D,
\end{equation}
where $\varphi(t)$ and $\psi(t)$ are nondecreasing functions of the
nonnegative argument. If $\varphi(t)$ and $\psi(t)$ are constants,
the conditions \eqref{u1} and \eqref{u2} are reduced to the classic
Lipschitz conditions.

In the case when $g = 0$ the most precise error estimates for the
process \eqref{priblNewton} were obtained in \cite{GW1, GW2} under
a new smoothness assumption imposed on the operator $f$ called
regular smoothness. In this parer we generalize the main result
from \cite{GW2} to equations of the form \eqref{osnur} under the
hypotheses that the operator $f$ is regularly smooth on $D$ and
the operator $g$ satisfies \eqref{u2}. The convergence theorem for
the process \eqref{priblNewton} is proved by means of majorant
equations.

\section{Regular smoothness}

Let $\mathcal{N}$ denote the class of continuous strictly increasing
functions $\omega: [0, \infty) \rightarrow [0, \infty)$ that are
concave and vanishing at zero: $\omega(0) = 0.$ Assume without loss
of generality that $f'(x_0) = I.$

Denote by $h(f)$ the quantity $\inf\limits_{x \in D}\|f'(x)\|.$
Given an $\omega \in \mathcal{N},$ we say in accordance with
\cite{GW2} that $f$ is \textit{$\omega$-regularly smooth} on $D$
(or, equivalently, that $\omega$ is a \textit{regular smoothness
modulus of $f$} on $D),$ if there exists $h \in [0, h(f)]$ such that
the inequality
\begin{equation}\label{neravregul1}
\omega^{-1}\left(h_f(x', x'') + \|f'(x'') - f'(x')\|\right) -
\omega^{-1}\left(h_f(x', x'')\right) \leq \|x'' - x'\|,
\end{equation}
where
$$
h_f(x', x'') = \min\{\|f'(x')\|,\ \|f'(x'')\|\} - h,
$$
holds for all $x',\, x'' \in D.$

The operator $f$ is called \textit{regularly smooth} on $D,$ if it
is $\omega$-regularly smooth on $D$ for some $\omega \in
\mathcal{N}.$

The condition \eqref{neravregul1} may be written in the form
$$
\|f'(x'') - f'(x')\| \leq \omega\left(\omega^{-1}(h_f(x', x'')) +
\|x'' - x'\|\right) - h_f(x', x''),
$$
or
\begin{equation}\label{neravregul2}
\|f'(x'') - f'(x')\| \leq \omega(\xi(x', x'') + \|x'' - x'\|) -
\omega(\xi(x', x'')),
\end{equation}
where $\xi(x', x'') = \omega^{-1}(h_f(x', x'')).$

It should be remarked that in \cite{GW1} a more restrictive
definition of regular smoothness was used, which coincides with the
definition in \cite{GW2} when $h = 0.$ In fact, if for some $h =
h_0$ and some $x',\ x'' \in D$ the inequality \eqref{neravregul2}
holds, then it will be true for all $h > h_0$ with the same $x',\
x'' \in D$ because of the difference $\omega(t + \tau) - \omega(t)$
does not increase in $t$ for each fixed $\tau
> 0.$

\begin{lemma}\label{lemma:regul} {\rm \cite{GW2}}
If the operator $f$ is $\omega$-regularly smooth on $D$ with some
$h,$ then
$$
\left|\omega^{-1}(\|f'(x'')\| - h) - \omega^{-1}(\|f'(x')\| -
h)\right| \leq \|x'' - x'\|
$$
for all $x',\ x'' \in D.$
\end{lemma}

It follows from the definition of $\xi$ and
Lemma~\ref{lemma:regul} that
\begin{equation}\label{neravksi}
\xi(x', x'') \geq \omega^{-1}(\|f'(x')\| - h) - \|x'' - x'\|
\end{equation}
for all $x',\ x'' \in D.$

\section{Some preliminary results}

The proof of the main theorem is based on several preliminary
propositions.

Let $\omega \in \mathcal{N},$ $\Omega(t) = \int\limits_0^t
\omega(\tau)\,d\tau,$ $\Psi(t) = \int\limits_0^t \psi(\tau)\,d\tau,$
$\chi = \omega^{-1}(1 - h),$ $a$ is a positive number such that
$$
\|f(x_0) + g(x_0)\| \leq a
$$
and
$$
\Phi_h(t) = a - \Omega(\chi) + \Omega(\chi - t) - th,\qquad t \in
[0,\, \chi].
$$

Let us define the numerical sequence $\{t_n\}$ by the following
recurrence formula:
\begin{equation}\label{posledov}
t_{n+1} = t_n + \D\frac{a - \Omega(\chi) + \Omega(\chi - t_n) -
t_nh + \Psi(t_n)}{h + \omega(\chi - t_n)},
\end{equation}
$n = 0,\ 1,\ \ldots;$ $t_0 = 0.$

In terms of the function
\begin{equation}\label{funcW}
W(t) = \Phi_h(t) + \Psi(t)
\end{equation}
the relation \eqref{posledov} may be rewritten as follows:
$$
t_{n+1} = t_n - \D\frac{W(t_n)}{\Phi'_h(t_n)},
$$
$n = 0,\ 1,\ \ldots;$ $t_0 = 0.$

\begin{lemma}\label{lemma:posledov}
Suppose that the function \eqref{funcW} has a unique zero $t_\ast$
in the interval $[0, \chi]$ and
\begin{equation}\label{uslW}
a < \Omega(\chi) + h \cdot \chi - \Psi(\chi).
\end{equation}
Then the sequence \eqref{posledov} is defined for all $n,$
monotonically increases and converges to $t_\ast.$
\end{lemma}

\dok The function $W$ is positive on the interval $[0, t_\ast),$
since $t_\ast$ is a unique zero of the equation $W(t) = 0,$ $W(0) =
a > 0$ and $W$ is continuous on $[0, \chi].$ Moreover, the function
$\Phi'_h(t) = -\omega(\chi - t) - h$ is negative on the interval
$[0, t_\ast).$ Hence the function
$$
u(t) = - \D\frac{W(t)}{\Phi'_h(t)}
$$
is positive on $[0, t_\ast).$

Let us show that the function $t + u(t)$ is nondecreasing on $[0,
t_\ast).$ In fact,
$$
(t + u(t))' = 1 + u'(t) = 1 + \left(\D\frac{\Phi_h(t) +
\Psi(t)}{\omega(\chi - t) + h}\right)' =
$$
$$
= 1 + \D\frac{(\Phi'_h(t) + \Psi'(t)) \cdot (\omega(\chi - t) + h)
+ (\Phi_h(t) + \Psi(t)) \cdot \omega '(\chi - t)}{(\omega(\chi -
t) + h)^2} =
$$
$$
= \D\frac{\Psi'(t) \cdot (\omega(\chi - t) + h) + (\Phi_h(t) +
\Psi(t)) \cdot \omega '(\chi - t)}{(\omega(\chi - t) + h)^2} \geq
0
$$
on $[0, t_\ast).$ This implies that the sequence $\{t_n\}$
monotonically increases and $t_{n+1} = t_n + u(t_n) \leq t_\ast +
u(t_\ast) = t_\ast$ for $t_n \leq t_\ast.$ Consequently, the
sequence $\{t_n\}$ converges to $t_{\ast\ast} \in [0, t_\ast]$ and
$t_{\ast\ast} = t_{\ast\ast} + u(t_{\ast\ast}),$ hence
$W(t_{\ast\ast}) = 0.$ Since $t_\ast$ is a unique zero of $W$ in the
interval $[0, \chi],$ it follows that $t_{\ast\ast} = t_\ast.$

The  sequence $\{t_n\}$ is defined for all $n.$ In fact, it is clear
from \eqref{uslW} that $W(\chi) < 0 < a = W(0)$ and hence there
exists $\theta \in (0, \chi)$ such that $W(\theta) = 0.$
Consequently, $\theta = t_\ast = \lim\limits_{n \to \infty}t_n$ and
$t_n \leq \theta < \chi$ for all $n = 0,\ 1,\ \ldots.$ Because of
the monotonicity of $\omega$ the inequality $\omega(\chi - t_n)
> 0$ is true for all $n = 0,\ 1,\ \ldots.$ This completes the proof
of Lemma~\ref{lemma:posledov}.

\begin{lemma}\label{lemma:reszur}
Let the operator $f$ be $\omega$-regularly smooth on $D$ with some
$h,$ the operator $g$ satisfies \eqref{u2}, the function
\eqref{funcW} has a unique zero $t_\ast$ in the interval $[0, \chi]$
and the closed ball $\overline{B(x_0, t_\ast)}$ is contained in $D.$
Then the equation \eqref{osnur} has a unique solution $x_\ast$ in
the ball $\overline{B(x_0, t_\ast)}.$
\end{lemma}

\dok Let us prove the existence of a solution in the ball
$\overline{B(x_0, t_\ast)}.$ Consider the sequence
$$
u_{n+1} = Du_n\qquad (n = 0,\ 1,\ \ldots;\ u_0 = x_0),
$$
where $D = I - [f'(x_0)]^{-1}(f + g) = I - (f + g),$ and the
numerical sequence
$$
\rho_{n+1} = d(\rho_n)\qquad (n = 0,\ 1,\ \ldots;\ \rho_0 = 0),
$$
where $d(t) = t + W(t).$ Since
$$
d'(t) = 1 + W'(t) = 1 + \Phi'_h(t) + \Psi'(t) =
$$
$$
= 1 - h - \omega(\chi - t) + \psi(t) = \omega(\chi) - \omega(\chi
- t) + \psi(t) \geq 0
$$
on the interval $[0, \chi],$ the function $d$ is monotonically
increasing on $[0, \chi].$

For all $n = 0,\ 1,\ \ldots$ the inequality
\begin{equation}\label{inrho}
\rho_n \leq t_\ast
\end{equation}
holds. In fact, for $n = 0$ the inequality \eqref{inrho} is obvious:
$\rho_0 = 0 \leq t_\ast.$ Suppose that \eqref{inrho} holds for all
$n \leq k.$ Then from $\rho_k \leq t_\ast$ because of the
monotonicity of $d$ we obtain $d(\rho_k) \leq d(t_\ast),$ that is
$\rho_{k + 1} \leq t_\ast.$ Consequently, by the induction
hypothesis the inequality \eqref{inrho} is true for all $n.$

Let us prove by induction that the sequence $\{\rho_n\}$ is
monotone. Clearly $0 = \rho_0 \leq \rho_1 = a.$ Suppose that $\rho_k
\leq \rho_{k + 1}.$ Then $\rho_{k + 1} = d(\rho_k) \leq d(\rho_{k +
1}) = \rho_{k + 2}.$

Thus the sequence $\{\rho_n\}$ is monotonically increasing and
bounded from above. Consequently, it converges to some $\tilde{\rho}
\in [0, t_\ast].$ By letting $n \rightarrow \infty$ in $\rho_{n + 1}
= \rho_n + W(\rho_n)$ we obtain $W(\tilde{\rho}) = 0$ and
$\tilde{\rho} = t_\ast.$

Let us show that for all $n = 0,\ 1,\ \ldots$ the inequality
\begin{equation}\label{nerrho}
\|u_{n+1} - u_n\| \leq \rho_{n+1} - \rho_n
\end{equation}
holds.

For $n = 0$ the inequality \eqref{nerrho} is obvious:
$$
\|u_1 - u_0\| = \|x_0 - (f(x_0) + g(x_0)) - x_0\| = \|f(x_0) +
g(x_0)\| \leq a = W(0) = \rho_1 - \rho_0.
$$
Suppose that \eqref{nerrho} holds for all $n < k.$ Then
$$
\|u_{k+1} - u_k\| = \|Du_k - Du_{k-1}\| = \|u_k - u_{k-1} -
(f(u_k) - f(u_{k-1})) - (g(u_k) - g(u_{k-1}))\| \leq
$$
$$
\leq \|u_k - u_{k-1} - (f(u_k) - f(u_{k-1}))\| + \|g(u_k) -
g(u_{k-1})\| \leq
$$
$$
\leq \int\limits_0^1 \|f'(u_t) - f'(x_0)\|\|u_k - u_{k-1}\|\,dt +
\|g(u_k) - g(u_{k - 1})\| \leq
$$
$$
\leq \int\limits_0^1 (\omega(\xi(x_0, u_t) + \|u_t - x_0\|) -
\omega(\xi(x_0, u_t)))\|u_k - u_{k-1}\|\,dt + \|g(u_k) - g(u_{k -
1})\|,
$$
where $u_t = u_{k-1} + t(u_k - u_{k-1}),$ $0 \leq t \leq 1.$

By the inequality \eqref{neravksi} we have
$$
\xi(x_0, u_t) \geq \omega^{-1}(\|f'(x_0)\| - h) - \|u_t - x_0\| =
\chi - \|u_t - x_0\|.
$$
By the induction hypothesis
$$
\|u_k - x_0\| = \|u_k - u_0\| \leq \sum_{j=1}^k \|u_j- u_{j-1}\|
\leq \sum_{j=1}^k (\rho_j - \rho_{j-1}) = \rho_k.
$$
Consequently,
$$
\|u_t - x_0\| = \|(1-t)(u_{k-1} - u_0) + t(u_k - u_0)\| \leq
(1-t)\|u_{k-1} - u_0\| + t\|u_k - u_0\| \leq
$$
$$
\leq (1-t)\rho_{k-1} + t\rho_k.
$$

From \eqref{u2} and Proposition~1 in \cite{ZF} it follows that
\begin{equation}\label{u3}
\|g(x'') - g(x')\| \leq \Psi(t + \|x'' - x'\|) - \Psi(t)\quad
\forall\,x',\ x'' \in \overline{B(x_0, t)} \subseteq D.
\end{equation}

Because of concavity of $\omega$ and \eqref{u3} we have
$$
\|u_{k+1} - u_k\| \leq \int\limits_0^1 (\omega(\chi - \|u_t -
x_0\| + \|u_t - x_0\|) - \omega(\chi - \|u_t - x_0\|))\|u_k -
u_{k-1}\|\,dt +
$$
$$
+ \Psi(\rho_{k-1} + \|u_k - u_{k-1}\|) - \Psi(\rho_{k-1}) \leq
$$
$$
\leq \int\limits_0^1 (\omega(\chi) - \omega(\chi - \|u_t -
x_0\|))(\rho_k - \rho_{k-1})\,dt + \Psi(\rho_k) - \Psi(\rho_{k-1})
\leq
$$
$$
\leq \int\limits_0^1 (\omega(\chi) - \omega(\chi -
((1-t)\rho_{k-1} + t\rho_k)))(\rho_k - \rho_{k-1})\,dt +
\Psi(\rho_k) - \Psi(\rho_{k-1}) =
$$
$$
= \int\limits_0^1 (1 + \Phi_h'((1-t)\rho_{k-1} + t\rho_k))(\rho_k
- \rho_{k-1})\,dt + \Psi(\rho_k) - \Psi(\rho_{k-1}) =
$$
$$
= \int\limits_{\rho_{k-1}}^{\rho_k}(1 + \Phi_h'(\theta))\,d\theta
+ \Psi(\rho_k) - \Psi(\rho_{k-1}) =
$$
$$
= \rho_k - \rho_{k-1} + \Phi_h(\rho_k) - \Phi_h(\rho_{k-1}) +
\Psi(\rho_k) - \Psi(\rho_{k-1}) = d(\rho_k) - d(\rho_{k-1}) =
\rho_{k+1} - \rho_k.
$$
Thus the inequality \eqref{nerrho} holds for $n = k.$

It follows from \eqref{nerrho} that for $m > n$
$$
\|u_m - u_n\| \leq \|u_m - u_{m-1}\| + \ldots + \|u_{n+1} - u_n\|
\leq \rho_m - \rho_{m-1} + \ldots + \rho_{n+1} - \rho_n = \rho_m -
\rho_n.
$$
Hence for all $m$ and $n$
\begin{equation}\label{equ}
\|u_m - u_n\| \leq |\rho_m - \rho_n|.
\end{equation}

Since the sequence $\{\rho_n\}$ converges to $t_\ast,$ it follows
from \eqref{equ} that the sequence $\{u_n\}$ also converges to some
$x_\ast.$ Further,
$$
\|u_n - u_0\| \leq \rho_n \leq t_\ast\qquad (n = 0,\ 1,\ \ldots)
$$
and, consequently, all $u_n$ with $x_\ast$ belong to the ball
$\overline{B(x_0, t_\ast)}.$ By letting $n \rightarrow \infty$ in
$u_{n+1} = Du_n$ we obtain that $x_\ast = D(x_\ast)$ or $f(x_\ast) +
g(x_\ast) = 0.$ Thus $x_\ast$ is a solution of the equation
\eqref{osnur} in the ball $\overline{B(x_0, t_\ast)}.$

To prove the uniqueness of the solution $x_\ast$ in the ball
$\overline{B(x_0, t_\ast)}$ consider the second solution
$x_{\ast\ast} \in \overline{B(x_0, t_\ast)}$ of \eqref{osnur} and
show that for all $n = 0,\ 1,\ \ldots$ the inequality
\begin{equation}\label{ner3}
\|x_{\ast\ast} - u_n\| \leq t_\ast - \rho_n
\end{equation}
holds.

For $n = 0$ the inequality \eqref{ner3} is obvious:
$$
\|x_{\ast\ast} - x_0\| \leq t_\ast - \rho_0 = t_\ast.
$$
Suppose that \eqref{ner3} holds for all $n \leq k.$ Then
$$
\|x_{\ast\ast} - u_{k+1}\| = \|x_{\ast\ast} - Du_k\| =
\|x_{\ast\ast} - u_k + f(u_k) + g(u_k)\| =
$$
$$
= \|f(u_k) - f(x_{\ast\ast}) - (u_k - x_{\ast\ast}) + g(u_k) -
g(x_{\ast\ast})\| \leq
$$
$$
\leq \|f(u_k) - f(x_{\ast\ast}) - f'(x_0)(u_k - x_{\ast\ast})\| +
\|g(u_k) - g(x_{\ast\ast})\| \leq
$$
$$
\leq \int\limits_0^1 \|f'(\tilde{u}_t) - f'(x_0)\|\|u_k -
x_{\ast\ast}\|\,dt + \|g(u_k) - g(x_{\ast\ast})\| \leq
$$
$$
\leq \int\limits_0^1 (\omega(\xi(x_0, \tilde{u}_t) + \|\tilde{u}_t
- x_0\|) - \omega(\xi(x_0, \tilde{u}_t)))\,\|u_k -
x_{\ast\ast}\|\,dt + \|g(u_k) - g(x_{\ast\ast})\|,
$$
where $\tilde{u}_t = x_{\ast\ast} + t(u_k - x_{\ast\ast}),$ $0 \leq
t \leq 1.$

By the inequality \eqref{neravksi} we have
$$
\xi(x_0, \tilde{u}_t) \geq \omega^{-1}(\|f'(x_0)\| - h) -
\|\tilde{u}_t - x_0\| = \chi - \|\tilde{u}_t - x_0\|.
$$

Further,
$$
\|\tilde{u}_t - x_0\| = \|(1-t)(x_{\ast\ast} - x_0) + t(u_k -
x_0)\| \leq (1-t)\|x_{\ast\ast} - x_0\| + t\|u_k - x_0\| \leq
$$
$$
\leq (1-t)t_\ast + t\rho_k.
$$
Because of concavity of $\omega,$ the inequality \eqref{u3} and the
induction hypothesis we have
$$
\|x_{\ast\ast} - u_{k+1}\| \leq
$$
$$
\leq \int\limits_0^1 (\omega(\chi - \|\tilde{u}_t - x_0\| +
\|\tilde{u}_t - x_0\|) - \omega(\chi - \|\tilde{u}_t -
x_0\|))\,\|u_k - x_{\ast\ast}\|\,dt + \|g(u_k) - g(x_{\ast\ast})\|
\leq
$$
$$
\leq \int\limits_0^1 (\omega(\chi) - \omega(\chi - \|\tilde{u}_t -
x_0\|))(t_\ast - \rho_k)\,dt + \Psi(\rho_k + \|u_k -
x_{\ast\ast}\|) - \Psi(\rho_k) \leq
$$
$$
\leq \int\limits_0^1 (\omega(\chi) - \omega(\chi - ((1-t)t_\ast +
t\rho_k)))(t_\ast - \rho_k)\,dt + \Psi(t_\ast) - \Psi(\rho_k) =
$$
$$
= \int\limits_0^1 (1 + \Phi'_h((1-t)t_\ast + t\rho_k))(t_\ast -
\rho_k)\,dt + \Psi(t_\ast) - \Psi(\rho_k) =
$$
$$
= \int\limits_{\rho_k}^{t_\ast}(1 + \Phi'_h(\theta))\,d\theta +
\Psi(t_\ast) - \Psi(\rho_k) =
$$
$$
= t_\ast - \rho_k + \Phi_h(t_\ast) - \Phi_h(\rho_k) + \Psi(t_\ast)
- \Psi(\rho_k) = d(t_\ast) - d(\rho_k) = t_\ast - \rho_{k+1}.
$$
Hence \eqref{ner3} holds for $n = k+1.$

By letting $n \rightarrow \infty$ in \eqref{ner3} we obtain that
$$
\|x_{\ast\ast} - x_\ast\| \leq t_\ast - t_\ast = 0
$$
and hence $x_{\ast\ast} = x_\ast.$ This completes the proof of
Lemma~\ref{lemma:reszur}.

Let us denote for all $n = 1,\ 2,\ \ldots$
$$
r(x_{n-1}, x_n) = \|f(x_n) - f(x_{n-1}) - f'(x_{n-1})(x_n -
x_{n-1})\|.
$$

\begin{lemma}\label{lemma:ocenka}
Let the operator $f$ be $\omega$-regularly smooth on $D$ with some
$h,$ the operator $g$ satisfies \eqref{u2}, the sequence $\{t_n\}$
is defined by the recurrence formula \eqref{posledov} and the
condition \eqref{uslW} holds. If for all $1 \leq k \leq n$
successive approximations $x_k$ are defined and satisfy the
inequality
\begin{equation}\label{nerav0}
\|x_k - x_{k-1}\| \leq t_k - t_{k-1},
\end{equation}
then
\begin{equation}\label{ocenka}
r(x_{n-1}, x_n) \leq a - \Omega(\chi) + \Omega(\chi - t_n) - t_nh
+ \Psi(t_{n-1}).
\end{equation}
\end{lemma}

\dok Let $x_t = x_{n-1} + t(x_n - x_{n-1}),$ $0 \leq t \leq 1.$ Then
$$
r(x_{n-1}, x_n) \leq \int\limits_0^1 \|f'(x_t) -
f'(x_{n-1})\|\|x_n - x_{n-1}\|\,dt \leq
$$
$$
\leq \int\limits_0^1 (\omega(\xi(x_{n-1}, x_t) + \|x_t -
x_{n-1}\|) - \omega(\xi(x_{n-1}, x_t)))\,\|x_n - x_{n-1}\|\,dt.
$$

By the inequality \eqref{neravksi} we have
$$
\xi(x_{n-1}, x_t) \geq \omega^{-1}(\|f'(x_{n-1})\| - h) - \|x_t -
x_{n-1}\| =
$$
$$
= \omega^{-1}(\|f'(x_{n-1})\| - h) - t\|x_n - x_{n-1}\|.
$$

Since for all $1 \leq k \leq n$ the inequality \eqref{nerav0}
holds, it follows that
$$
\|x_n - x_0\| \leq \sum_{k=1}^n \|x_k - x_{k-1}\| \leq
\sum_{k=1}^n (t_k - t_{k-1}) = t_n.
$$

By Lemma~\ref{lemma:regul}
$$
\omega^{-1}(\|f'(x_n)\| - h) \geq \omega^{-1}(\|f'(x_0)\| - h) -
\|x_n - x_0\|.
$$
Since $\omega^{-1}(\|f'(x_n)\| - h) \geq 0,$ we have
\begin{multline}\label{neravomega0}
\omega^{-1}(\|f'(x_n)\| - h) \geq \left(\omega^{-1}(\|f'(x_0)\| -
h) - \|x_n - x_0\|\right)^{+} \geq \\
\geq \left(\omega^{-1}(\|f'(x_0)\| - h) - t_n\right)^{+},
\end{multline}
where $\lambda^{+} = \max\{\lambda, 0\}.$

Further $\left(\omega^{-1}(\|f'(x_0)\| - h) - t_n\right)^{+} =
\left(\omega^{-1}(1 - h) - t_n\right)^{+} = (\chi - t_n)^{+}.$ Let
$\alpha_n = (\chi - t_n)^{+}.$ By the condition \eqref{uslW} we
obtain that $t_n < \chi$  and $\alpha_n = \chi - t_n > 0$ for all $n
= 0,\ 1,\ \ldots.$ Hence the inequality \eqref{neravomega0} may be
rewritten in the form
$$
\omega^{-1}(\|f'(x_n)\| - h) \geq \alpha_n.
$$
Analogously we obtain
$$
\omega^{-1}(\|f'(x_{n-1})\| - h) \geq \alpha_{n-1}.
$$

Using $\|x_n - x_{n-1}\| \leq t_n - t_{n-1}$ we get
$$
\omega^{-1}(\|f'(x_{n-1})\| - h) - t\|x_n - x_{n-1}\| \geq
\alpha_{n-1} - t(t_n - t_{n-1}),
$$
which implies that
$$
\xi(x_{n-1}, x_t) \geq \alpha_{n-1} - t(t_n - t_{n-1}) =
\alpha_{n-1} - t\delta_{n-1},
$$
where $\delta_{n-1} = t_n - t_{n-1}.$

Because of concavity and monotonicity of $\omega$ we have
$$
r(x_{n-1}, x_n) \leq \int\limits_0^1 (\omega(\alpha_{n-1} -
t\delta_{n-1} + t\delta_{n-1}) - \omega(\alpha_{n-1} -
t\delta_{n-1}))\delta_{n-1}\,dt =
$$
$$
= \int\limits_0^1 (\omega(\alpha_{n-1}) - \omega(\alpha_{n-1} -
t\delta_{n-1}))\delta_{n-1}\,dt.
$$

Let $\tau = t\delta_{n-1}.$ Then
$$
r(x_{n-1}, x_n) \leq
\int\limits_0^{\delta_{n-1}}(\omega(\alpha_{n-1}) -
\omega(\alpha_{n-1} - \tau))\,d\tau =
\int\limits_0^{\delta_{n-1}}\omega(\alpha_{n-1})\,d\tau -
\int\limits_0^{\delta_{n-1}}\omega(\alpha_{n-1} - \tau)\,d\tau =
$$
$$
= \omega(\alpha_{n-1})\delta_{n-1} +
\int\limits_{\alpha_{n-1}}^{\alpha_{n-1} -
\delta_{n-1}}\omega(\theta)\,d\theta =
\omega(\alpha_{n-1})\delta_{n-1} + \int\limits_0^{\alpha_{n-1} -
\delta_{n-1}}\omega(\theta)\,d\theta -
\int\limits_0^{\alpha_{n-1}}\omega(\theta)\,d\theta =
$$
$$
= \omega(\alpha_{n-1})\delta_{n-1} + \Omega(\alpha_{n-1} -
\delta_{n-1}) - \Omega(\alpha_{n-1}) =
\omega(\alpha_{n-1})\delta_{n-1} + \Omega(\alpha_n) -
\Omega(\alpha_{n-1}) =
$$
$$
= \omega(\chi - t_{n-1}) \cdot (t_n - t_{n-1}) + \Omega(\chi -
t_n) - \Omega(\chi - t_{n-1}).
$$

Let us show that for all $n = 0,\ 1,\ \ldots$ the equality
\begin{equation}\label{ravenstvo}
\omega(\chi - t_n) \cdot (t_{n+1} - t_n) - \Omega(\chi - t_n) +
t_{n+1} h - \Psi(t_n) = a - \Omega(\chi)
\end{equation}
holds. In fact, by the definition of the sequence $\{t_n\}$
$$
(t_{n+1} - t_n)(h + \omega(\chi - t_n)) = a - \Omega(\chi) +
\Omega(\chi - t_n) - t_nh + \Psi(t_n)
$$
and
$$
(t_n - t_{n-1})(h + \omega(\chi - t_{n-1})) = a - \Omega(\chi) +
\Omega(\chi - t_{n-1}) - t_{n-1}h + \Psi(t_{n-1}).
$$
It follows from the first of these equalities that
$$
a - \Omega(\chi) = t_{n+1}h + \omega(\chi - t_n) \cdot (t_{n+1} -
t_n) - \Omega(\chi - t_n) - \Psi(t_n)
$$
and from the second that
$$
a - \Omega(\chi) = t_nh + \omega(\chi - t_{n-1}) \cdot (t_n -
t_{n-1}) - \Omega(\chi - t_{n-1}) - \Psi(t_{n-1}).
$$
Consequently,
$$
\omega(\chi - t_n) \cdot (t_{n+1} - t_n) - \Omega(\chi - t_n) +
t_{n+1}h - \Psi(t_n) =
$$
$$
= \omega(\chi - t_{n-1}) \cdot (t_n - t_{n-1}) - \Omega(\chi -
t_{n-1}) + t_nh - \Psi(t_{n-1})
$$
for all $n = 1,\ 2,\ \ldots$ and
$$
\omega(\chi - t_n) \cdot (t_{n+1} - t_n) - \Omega(\chi - t_n) +
t_{n+1}h - \Psi(t_n) =
$$
$$
= \omega(\chi - t_0) \cdot (t_1 - t_0) - \Omega(\chi - t_0) + t_1h
- \Psi(t_0) =
$$
$$
= \omega(\chi) \cdot a - \Omega(\chi) + ah = (1 - h)a -
\Omega(\chi) + ah = a - \Omega(\chi).
$$
Thus the equality \eqref{ravenstvo} holds for all $n = 0,\ 1,\
\ldots$ and the estimate for $r(x_{n-1}, x_n)$ may be rewritten in
the form \eqref{ocenka}. This completes the proof of
Lemma~\ref{lemma:ocenka}.

\section{Convergence Theorem}

{\it Let the operator $f$ be $\omega$-regularly smooth on $D$ with
some $h,$ the operator $g$ satisfies \eqref{u2}, the function
\eqref{funcW} has a unique zero $t_\ast$ in the interval $[0,
\chi],$ the closed ball $\overline{B(x_0, t_\ast)}$ is contained in
$D$ and the condition \eqref{uslW} holds. Then

\begin{enumerate}[\hskip 0.6 cm \rm 1)]

\item the equation \eqref{osnur} has a unique root $x_\ast$ in the
ball $\overline{B(x_0, t_\ast)};$

\item the successive approximations \eqref{priblNewton} are defined
for all $n = 0,\ 1,\ \ldots,$ belong to $\overline{B(x_0, t_\ast)}$
and converge to $x_\ast;$

\item for all $n = 0,\ 1,\ \ldots$ the inequalities
\begin{equation}\label{nerav1}
\|x_{n+1} - x_n\| \leq t_{n+1} - t_n,
\end{equation}
\begin{equation}\label{nerav2}
\|x_\ast - x_n\| \leq t_\ast - t_n,
\end{equation}
hold, where the sequence $\{t_n\}$ is defined by the recurrence
formula \eqref{posledov}, monotonically increases and converges to
$t_\ast.$
\end{enumerate}
}

\dok In order to prove the theorem it suffices to show that the
successive approximations \eqref{priblNewton} are defined for all $n
= 0,\ 1,\ \ldots,$ belong to the ball $\overline{B(x_0, t_\ast)}$
and satisfy the inequalities \eqref{nerav1} and \eqref{nerav2}.
Other assertions of the theorem follow from
Lemma~\ref{lemma:posledov} and Lemma~\ref{lemma:reszur}.

Since \eqref{nerav2} is a direct consequence of \eqref{nerav1}, it
suffices to prove \eqref{nerav1}. For $n = 0$ the inequality
\eqref{nerav1} is obvious:
$$
\|x_1 - x_0\| = \left\|[f'(x_0)]^{-1}(f(x_0) + g(x_0))\right\|
\leq a = t_1 - t_0.
$$

Suppose that \eqref{nerav1} holds for all $n < k.$ We first show
that the operator $f'(x_k)$ is invertible. In fact,
$$
\left\|[f'(x_0)]^{-1}(f'(x_k) - f'(x_0))\right\| = \|f'(x_k) -
f'(x_0)\| \leq
$$
$$
\leq \omega(\xi(x_0, x_k) + \|x_k - x_0\|) - \omega(\xi(x_0,
x_k)).
$$
By the inequality \eqref{neravksi}
$$
\xi(x_0, x_k) \geq \omega^{-1}(\|f'(x_0)\| - h) - \|x_k - x_0\| =
\chi - \|x_k - x_0\|.
$$
By the induction hypothesis
$$
\|x_k - x_0\| \leq \sum_{j=1}^k \|x_j- x_{j-1}\| \leq \sum_{j=1}^k
(t_j - t_{j-1}) = t_k
$$
and hence $\xi(x_0, x_k) \geq \chi - t_k > 0$ ($t_k < \chi$ for
all $k = 0,\ 1,\ \ldots$ due to \eqref{uslW}). Because of
concavity of $\omega$ we have
$$
\omega(\xi(x_0, x_k) + \|x_k - x_0\|) - \omega(\xi(x_0, x_k)) \leq
\omega(\chi - t_k + \|x_k - x_0\|) - \omega(\chi - t_k) \leq
$$
$$
\leq \omega(\chi - t_k + t_k) - \omega(\chi - t_k) < \omega(\chi)
- \omega(0) = \omega(\chi) = 1 - h \leq 1.
$$

Thus $\left\|[f'(x_0)]^{-1}(f'(x_k) - f'(x_0))\right\| < 1$ and,
consequently, the operator
$$
T = I + [f'(x_0)]^{-1}(f'(x_k) - f'(x_0))
$$
is invertible. Since $f'(x_k) = f'(x_0)T = T,$ the operator
$f'(x_k)$ is also invertible and
$$
\left\|[f'(x_k)]^{-1}\right\| = \left\|T^{-1}\right\| \leq
\frac{1}{1 - \|T - I\|} \leq \frac{1}{1 - [\omega(\chi) -
\omega(\chi - t_k)]}.
$$

Further, using the estimate for $r(x_{k-1}, x_k)$ from
Lemma~\ref{lemma:ocenka} and the inequality \eqref{u3} we get
$$
\|x_{k+1} - x_k\| = \left\|[f'(x_k)]^{-1}(f(x_k) + g(x_k))\right\|
=
$$
$$
= \left\|[f'(x_k)]^{-1}(f(x_k) - f(x_{k-1}) - f'(x_{k-1})(x_k -
x_{k-1}) + g(x_k) - g(x_{k-1}))\right\| \leq
$$
\begin{multline*}
\leq \left\|[f'(x_k)]^{-1}\right\| \cdot \|f(x_k) - f(x_{k-1}) -
f'(x_{k-1})(x_k - x_{k-1})\| + \\ + \left\|[f'(x_k)]^{-1}\right\|
\cdot \|g(x_k) - g(x_{k-1})\| \leq
\end{multline*}
$$
\leq \D\frac{r(x_{k-1}, x_k) + \Psi(t_{k-1} + \|x_k - x_{k-1}\|) -
\Psi(t_{k-1})}{1 - [\omega(\chi) - \omega(\chi - t_k)]} \leq
$$
$$
\leq \D\frac{r(x_{k-1}, x_k) + \Psi(t_k) - \Psi(t_{k-1})}{1 -
[\omega(\chi) - \omega(\chi - t_k)]} \leq
$$
$$
\leq \D\frac{a - \Omega(\chi) + \Omega(\chi - t_k) - t_kh +
\Psi(t_{k-1}) + \Psi(t_k) - \Psi(t_{k-1})}{h + \omega(\chi - t_k)}
=
$$
$$
= \D\frac{a - \Omega(\chi) + \Omega(\chi - t_k) - t_kh +
\Psi(t_k)}{h + \omega(\chi - t_k)} = t_{k+1} - t_k.
$$
Consequently, \eqref{nerav1} holds for $n = k.$

Since for all $n = 0,\ 1,\ \ldots$ the operator $f'(x_n)$ is
invertible and $\|x_n - x_0\| \leq t_n \leq t_\ast,$ the successive
approximations \eqref{priblNewton} are defined for all $n = 0,\ 1,\
\ldots$ and belong to the ball $\overline{B(x_0, t_\ast)}.$ The
convergence of successive approximations to $x_\ast$ follows from
\eqref{nerav2}. This proves the theorem.

It is to be noted that each Lipschitz smooth operator is also
regularly smooth, but the converse is not true. So the theorem
proved is applicable to more wide class of nonlinear operator
equations of the form \eqref{osnur} than the corresponding
convergence theorems from \cite{ZabrZlep, ZF}.

\section*{Acknowledgments}

The author is grateful to Professor Petr~P. Zabreiko for his careful
reading of the paper and valuable help in improving its contents.

\renewcommand{\refname}{References}

\end{document}